\documentclass[12pt]{article}
\usepackage{amsmath,amssymb,amsthm, amsfonts, bm}
\usepackage[mathscr]{eucal}
\usepackage{stackengine}
\usepackage{accents}
\usepackage{rotate,graphics,epsfig}
\usepackage{color}

\newtheorem{The}{Theorem}

\newtheorem{Pro}[The]{Proposition}

\theoremstyle{definition}
\newtheorem{Def}[The]{Definition}
\newtheorem{Rem}[The]{Remark}
\numberwithin{equation}{section}
\numberwithin{The}{section}
\numberwithin{figure}{section}

\newcommand{\be}{\begin{eqnarray}}
	\newcommand{\ee}{\end{eqnarray}}

\newcommand{\by}{\begin{eqnarray*}}
	\newcommand{\ey}{\end{eqnarray*}}

\newcommand{\bn}{\begin{enumerate}}
	\newcommand{\en}{\end{enumerate}}

\newcommand{\bi}{\begin{itemize}}
	\newcommand{\ei}{\end{itemize}}

\def\Fbar{{\overline F}}
\def\RV{\mbox{\rm{RV}}}

\usepackage[utf8]{inputenc}
\usepackage[T1]{fontenc}

\usepackage{color}

\begin{document}
\title{Operator Tail Densities of Multivariate Copulas}
\author{
Haijun Li
\footnote{{\small\texttt{lih@math.wsu.edu}}, Department of Mathematics and Statistics, Washington State University, 
	Pullman, WA 99164, U.S.A.}
}
\date{December 2025}
\maketitle

\begin{abstract}
	
Operator regular variation of a multivariate distribution can be decomposed into the operator tail dependence of the underlying copula and the regular variation of the univariate marginals. In this paper, we introduce operator tail densities for copulas and show that an operator-regularly-varying density can be characterized through the operator tail density of its copula together with the marginal regular variation. As an example, we demonstrate that although a Liouville copula is not available in closed form, it nevertheless admits an explicit operator tail-dependence function.
	
\end{abstract}

\section{Introduction}
\label{Instroduction}

Multivariate regular variation with tail-equivalent margins has been well studied (see, e.g., \cite{Resnick07}), and its structure can be decomposed into the tail dependence of the underlying copula together with the regular variation of the univariate marginals \cite{LS2009, Li2018}. When the univariate marginal distributions are not tail-equivalent, the copula approach, being univariate margin-free, can still be used to effectively extract scale-invariant tail-dependence information, which, combined with the univariate regular variation, facilitates the analysis of general multivariate regular variation.

Let $X=(X_1, \dots, X_d)$ be a non-negative random vector with the joint distribution $F$ and univariate marginal distributions $F_i$, $1\le i\le d$. The distribution $F$ (or $X$) is said to be multivariate regular varying with the intensity measure $\mu(\cdot)$, if there exists a univariate regularly varying function $U(t)$, as $t\to \infty$,  
\begin{equation}\label{MRV-e-1}
	\frac{\mathbb{P}(X\in tB)}{U(t)}\to \mu(B)
\end{equation}
for all relatively compact subset $B\subset \mathbb{R}_+^d\backslash\{0\}$ with $\mu(\partial B)=0$, bounded away from $0$. Due to the common scaling $t$ among all the components $X_i$'s used in \eqref{MRV-e-1}, the standard form \eqref{MRV-e-1} of multivariate regular variation is usually used under the assumption that the univariate marginal tails are equivalent in the sense of $(1-F_i(t))/(1-F_1(t))\to c_i>0$, for some constants $c_i$, $1\le i\le d$. For non-standard multivariate regular variation without tail-equivalent marginals, marginal transformations are typically applied to convert it into a standard form \cite{Resnick07}, and the most commonly used tool for this purpose is the copula method. Alternatively, one can replace the common scalar scaling $t$ in \eqref{MRV-e-1} with an operator scaling, which accommodates possibly different tail rates of the individual margins of a regularly varying random vector $X$.

It is worth noting that using different marginal tail scalings can reveal natural tail-dependence structures that are not visible under standard tail-equivalent scalings. In applications such as those involving financial data, multivariate heavy-tailed data often do not share a common tail index across all components, making it necessary to consider coordinate systems with operator norming to detect variations in tail behavior \cite{NPM01, MS99}. Li \cite{Li2018} introduced operator tail dependence for copulas to describe scale-invariant extremal dependence in multivariate distributions and showed that, for diagonal norming matrices, operator regular variation for a multivariate distribution and operator tail dependence for the underlying copula are equivalent. Moreover, employing different marginal tail scalings can uncover more natural hidden tail-dependence structures, as demonstrated in the case of Marshall–Olkin copulas \cite{Li2018}.

Most multivariate distributions are specified through their density functions, and therefore operator regular variation must be described in terms of the localized operator tail dependence of the associated copula. The goal of this paper is to introduce the operator tail density of a copula and to show that any operator regularly varying density on 
$\mathbb{R}_+^d$ 
can be decomposed into the operator tail density of the underlying copula together with the univariate marginal regularly varying densities. Note that operator regular variation with diagonal norming matrices reveals extremal dependence structures of multivariate distributions without requiring marginal tail equivalence, whereas the operator tail densities for  copulas capture scale-invariant tail dependence patterns that arise across a broader class of distributions, including multivariate rapidly varying distributions and others that may not be amenable to an operator-regular-variation description.

The rest of the paper is organized as follows. In Section 2, we introduce operator tail densities of copulas with respect to a matrix of tail indices. Section 3 presents the main results on the equivalence between operator regular variation and the operator tail dependence of their underlying copulas in terms of densities. Section 4 concludes with an example involving Liouville distributions and their copulas, along with additional remarks.

 In what follows, 
two functions $f,g: {\mathbb{R}}\to \mathbb{R}$ are said to be tail equivalent, denoted by  $f(x)\sim g(x)$ as $x\to a$, $a\in \overline{\mathbb{R}}=\mathbb{R}\cup \{+\infty\}$, if $\lim_{x\to a}[f(x)/g(x)]= 1$. A univariate Borel-measurable function $V: \mathbb{R}_+\to \mathbb{R}_+$ is said to be regularly varying at $\infty$ with tail index $\rho\in \mathbb{R}$, denoted by $V\in \mbox{RV}_\rho$, if $V(tx)/V(t)\to x^\rho$ as $t\to \infty$ for any $x>0$. A function $V: \mathbb{R}_+\to \mathbb{R}_+$ is said to be regularly varying at $0$ with tail index $\rho\in \mathbb{R}$, denoted by $V\in \mbox{RV}_\rho(0)$, if $V(t^{-1})\in \mbox{RV}_\rho$. 
In particular, a function $V\in \mbox{RV}_0$ (or $V\in \mbox{RV}_0(0)$) is said to be slowly varying at $\infty$ (or at $0$). See, for example,  \cite{BGT1987, Resnick07} for details on these notions and on the theory of regular variation.

\section{Operator Tail Densities of Copulas}
\label{Operator RV section}

In this section, we introduce the operator tail density for copulas admitting a density and show that operator regular variation of a multivariate density implies operator regular variation of the associated distribution.

Let $X=(X_1, \dots, X_d)$ be a non-negative random vector, having the distribution $F$ with continuous  univariate marginal distributions $F_i$, $1\le i\le d$. The copula $C$ on the support $[0,1]^d$ is defined as the distribution of $(F_1(X_1), \dots, F_d(X_d))$; that is, 
$C(u_1, \dots, u_d)=F(F_1^{-1}(u_1), \dots, F_d^{-1}(u_d))$, $(u_1, \dots, u_d)\in [0,1]^d$. The survival copula $\widehat{C}$ on the support $[0,1]^d$ is defined as the distribution of $(1-F_1(X_1), \dots, 1-F_d(X_d))$. Obviously the lower-orthant limits towards $(0, \dots, 0)$ for $C$ is the same as the upper-orthant limits towards $(1, \dots, 1)$ for $\widehat{C}$. 
See, for example, \cite{Joe97} for detailed discussions on copulas and their applications.

\begin{Def}\label{operator-d-1} Let $C: [0,1]^d\to [0,1]$ be a copula with density $c(\cdot)$. 
	\begin{enumerate}
		\item The upper tail density of $C$ with tail orders ${\kappa} = (\lambda_1, \dots, \lambda_d)>0$ is defined as the limit, denoted by  $\lambda_{C}(w; \kappa)$, such that for every $w=(w_1, \dots, w_d)\in [0,\infty)^d\backslash\{0\}$, 
		\begin{equation}\label{operator-e-1}
			\lambda_{C}(w; \kappa):=\lim_{u\to 0}\frac{c\big(1-r_1(u)w_1, \dots, 1-r_d(u)w_d\big)}{u^{1-\sum_{i=1}^d\lambda_i}\ell(u)}, 
		\end{equation}
	holds locally uniformly in $w$ for some functions $r_i\in \mbox{RV}_{\lambda_i}(0)$, where $\lambda_i>0$, $1\le i\le d$, and for some function $\ell(\cdot)$ that is slowly varying at $0$. 
	\item The lower tail density of $C$ with tail orders ${\kappa} = (\lambda_1, \dots, \lambda_d)$ is defined as the upper tail density $\lambda_{\widehat C}(w; \kappa)$ of the survival copula $\widehat C$. That is, for every $w=(w_1, \dots, w_d)\in [0,\infty)^d\backslash\{0\}$, 
	\begin{equation}\label{operator-e-2}
		\lambda_{\widehat C}(w; \kappa)=\lim_{u\to 0}\frac{c\big(r_1(u)w_1, \dots, r_d(u)w_d\big)}{u^{1-\sum_{i=1}^d\lambda_i}\ell(u)}, 
	\end{equation}
holds locally uniformly in $w$ for some functions $r_i\in \mbox{RV}_{\lambda_i}(0)$, where $\lambda_i>0$, $1\le i\le d$,  and for some function $\ell(\cdot)$ that is slowly varying at $0$. 
		\end{enumerate}
\end{Def}

When $\lambda_i = 1$, $1\le i\le d$, the tail density $\lambda_{C}(w; \kappa)$ reduces to the tail density function introduced in \cite{LW2013}. When $\lambda_1=\lambda_2= \cdots = \lambda_d$, the tail density $\lambda_{C}(w; \kappa)$ reduces to the tail density function with tail order $1/\lambda_i$, introduced in \cite{LH14}. The tail densities of of skew-elliptical distributions are analyzed in \cite{JL19, Li2021}.

The tail density enjoys the following property of quasihomogeneity.

\begin{Pro}\label{operator-p-1} Let $C: [0,1]^d\to [0,1]$ be a copula with  density $c(\cdot)$. If $C$ admits the tail density $\lambda_{C}(w; \kappa)$, where $\kappa = (\lambda_1, \dots, \lambda_d)$, for $\lambda_i>0$, $1\le i\le d$, then for any fixed $t>0$, 
\begin{equation}\label{operator-e-3}
\lambda_C(t^{\lambda_1}w_1, \dots, t^{\lambda_d}w_d;\kappa)= t^{1-\sum_{i=1}^d\lambda_i}\lambda_C(w_1, \dots, w_d; \kappa), 
\end{equation}	
for every $w=(w_1, \dots, w_d)\in [0,\infty)^d\backslash\{0\}$. 
\end{Pro}

\noindent
{\sl Proof.}
It follows from the local uniform convergence that
\begin{align*}
\lambda_C(t^{\lambda_1}w_1, \dots, t^{\lambda_d}w_d;\kappa) &= \lim_{u\to 0}\frac{c\big(1-r_1(u)t^{\lambda_1}w_1, \dots, 1-r_d(u)t^{\lambda_d}w_d\big)}{u^{1-\sum_{i=1}^d\lambda_i}\ell(u)} \\
&= t^{1-\sum_{i=1}^d\lambda_i}\lim_{ut\to 0}\frac{c\big(1-r_1(ut)w_1, \dots, 1-r_d(ut)w_d\big)}{(ut)^{1-\sum_{i=1}^d\lambda_i}\ell(ut)} \\
&= t^{1-\sum_{i=1}^d\lambda_i}\lambda_C(w_1, \dots, w_d; \kappa), 
\end{align*}
where $\ell(ut)\sim \ell(u)$ as $u\to 0$. 
\hfill $\Box$

To express the quasihomogeneity \eqref{operator-e-3} in terms of group invariance, we use operator scaling \cite{Li2023}. Given a $d\times d$ matrix $E$, we define the exponential matrix
\[\exp(E)=\sum_{k=0}^\infty \frac{E^k}{k!}, \ \mbox{where}\ E^0=I\ \mbox{(the $d\times d$ identity matrix)}, 
\]
and the power matrix
\begin{equation}
	\label{power matrix}
	t^E=\exp(E\log t)= \sum_{k=0}^\infty \frac{E^k(\log t)^k}{k!}, \ \mbox{for}\ t>0.
\end{equation}
Power matrices can be viewed as linear operators from $\mathbb{R}^d$ to $\mathbb{R}^d$, and behave like power functions; for example, $t^{-E}=(t^{-1})^E=(t^E)^{-1}$. For any positive-definite matrix $E$ and any norm $||\cdot||$ on $\mathbb{R}^d$, $||t^Ew||\to \infty$, as $t\to \infty$, uniformly on compact subsets of $w\in \mathbb{R}^d\backslash \{0\}$. Let $\mbox{GL}(\mathbb{R}^d)$ denote the Lie group of invertible $d\times d$ real matrices. Given a $d\times d$ real matrix $E$, a Borel measurable mapping $f: \mathbb{R}_+\to \mbox{GL}(\mathbb{R}^d)$ is said to be operator-regularly varying at infinity with index $E$ if 
\begin{equation}\label{operator RV}
	\lim_{t\to \infty}f(x t)f(t)^{-1} = x^E
\end{equation}
for all $x>0$. If $E=0$, then $f(t)$ is said to be operator-slowly varying at infinity. As in the case of multivariate regular variation, \eqref{operator RV} holds locally uniformly in $x>0$, which yields various representations for $f(\cdot)$. See \cite{MS01} for details on operator regular variation and its applications.

One drawback of the convergence \eqref{operator RV} is that the representation that $f(t) = t^EL(t)$, where $L(\cdot)$ is operator-slowly varying, does not hold for some $d\times d$ matrices $E$ \cite{MS01}. To avoid this technical issue, we focus on diagonal matrices $f(t)$. The convergence 
\eqref{operator RV} holds for a diagonal matrix $f(t)$ if and only if all the diagonal entries of $f(t)$ are univariate regularly varying and the index matrix $E$ is diagonal. Assume throughout this paper that the matrix $E$ is a diagonal matrix with real eigenvalues $\lambda_1, \dots, \lambda_d$: 
\begin{equation}
	E = 
	\begin{pmatrix}
		\lambda_1& \cdots & 0\\
		\vdots & \ddots & \vdots \\
		0 & \cdots & \lambda_d
	\end{pmatrix}. 
	\label{diag}
\end{equation}
It follows immediately that 
\begin{equation}
	t^E=
	\begin{pmatrix}
		t^{\lambda_1} & \cdots & 0\\
		\vdots & \ddots & \vdots \\
		0 & \cdots & t^{\lambda_d}
	\end{pmatrix},
	\ \ \mbox{for}\ t>0. 
	\label{dec}
\end{equation}
Therefore, the convergence 
\eqref{operator RV} holds for a diagonal matrix $f(t)$ if and only if 
\begin{equation}\label{operator RV-1}
	f(t) = t^EL(t),\ \mbox{and}\ 
		L(t)= \mbox{DIAG}(\ell_i(t)) := 
	\begin{pmatrix}
		\ell_1(t) & \cdots & 0\\
		\vdots & \ddots & \vdots \\
		0 & \cdots & \ell_d(t)
	\end{pmatrix},
	\ \ \mbox{for}\ t>0, 
\end{equation} 
where $\ell_i(t), 1\le i\le d$, are slowly varying at infinity. Here and hereafter $\mbox{DIAG}(a_i)$ denotes a diagonal matrix with the main diagonal entries $a_i$, $1\le i\le d$.

The quasihomogeneity \eqref{operator-e-3} can be rephrased as
\begin{equation}\label{operator-e-4}
\lambda_C(t^Ew; \kappa) = t^{1-\text{tr}(E)}\lambda_C(w; \kappa)
\end{equation}
  where $t^Ew = (t^{\lambda_1}w_1, \dots, t^{\lambda_d}w_d)$, for $w=(w_1, \dots, w_d)\in [0,\infty)^d\backslash\{0\}$. Here and hereafter $\text{tr}(E)$ denotes the trace of a matrix $E$. Observe that $(t^E, t>0)$ and $(t^{1-\text{tr}(E)}, t>0)$ are groups for a fixed $E$, and therefore, the quasihomogeneity \eqref{operator-e-3} entails the group invariance, that is fundamental in characterizing the tail density $\lambda_C(\cdot; \kappa)$.

Using power matrices as scaling, operator regular variation of a multivariate density is defined as follows.

\begin{Def}
	\label{d2.1} Suppose that a non-negative random vector $(X_1, \dots, X_d)$ has a distribution $F$ on $\mathbb{R}^d$ with density $f$. The density $f$ is said to be regularly varying with operator tail index $E$ given by \eqref{diag}, $\lambda_i>0$, $0\le i\le d$, denoted as $f\in \mbox{MRV}(E, -\rho,\lambda(\cdot))$, if there exists a mapping $g: \mathbb{R}_+\to \mbox{GL}(\mathbb{R}^d)$ that is operator-regularly varying, such that the convergence
	\begin{equation}
		\label{e2.1}
		\frac{f(g(t)x)}{t^{-\text{tr}(E)}V(t)}\to \lambda(x)>0, 
	\end{equation}
	holds locally uniformly in $x\in \mathbb{R}^d_+\backslash\{0\}$, for $V(t)\in \mbox{RV}_{-\rho}, \rho>0$. 
\end{Def}

 A comprehensive discussion on properties of exponential and power matrices, as well as operator regular variation in general can be found in \cite{MS01}. Operator regular variation for copulas has been studied in \cite{Li2018}.

\begin{Rem}\label{operator-r-1}
	\begin{enumerate}
		\item 		
	Operator regular variation has been defined and studied for a general matrix 
	$E$ in \cite{MS01}, where the eigenvalues of 
	$E$ were shown to play a central role. In this work, we focus on diagonal matrices 
	$E$ to more clearly illustrate how eigenvalues govern asymmetric joint tail behavior.
		\item The limiting function $\lambda(\cdot)$ satisfies the group invariance $\lambda(t^Ex) = t^{-\rho-\text{tr}(E)}\lambda(x)$, $x\in \mathbb{R}^d\backslash \{0\}$ and $t>0$. 
	\end{enumerate}
\end{Rem}

Clearly that $f\in \mbox{MRV}(E, -\rho,\lambda(\cdot))$, $\rho>0$, if and only if $f\in \mbox{MRV}(E/\rho, -1,\lambda(\cdot))$. Without loss of generality, we often discuss only the case that $f\in \mbox{MRV}(E, -1,\lambda(\cdot))$.

\begin{The}
	\label{pdf->CDF} 
	Let $X=(X_1, \dots, X_d)$ have a distribution $F$ defined on $\mathbb{R}^d$. If $F$ has a density $f\in \mbox{MRV}(E, -\rho,\lambda(\cdot))$ on $\{x: ||x||\ge \epsilon\}$, $\epsilon>0$, for any norm $||\cdot||$ on $\mathbb{R}^d$, where $\lambda(x)$ is locally bounded and  $E$ is given by \eqref{diag}, $\lambda_i>0$, $1\le i\le d$, then for any Borel subset $B\subseteq \{x: ||x||\ge \epsilon\}$, $\epsilon>0$,
	\begin{equation}
		\label{e2.3}
		\frac{\mathbb{P}(X\in t^EL(t)B)}{U(t)}\to \int_B\lambda(x)dx,
	\end{equation}
	where $L(t)$ is a diagonal matrix described in \eqref{operator RV-1}, and $U(t)\in \mbox{RV}_{-\rho}$.
\end{The} 

\noindent
{\sl Proof.} Since the density $f\in \mbox{MRV}(E, -\rho,\lambda(\cdot))$, 
there exists some $V(t)\in \mbox{RV}_{-\rho}$ and matrix $L(t)= \mbox{DIAG}(\ell_i(t))$, $\ell_i(t)\in \mbox{RV}_0$, $1\le i\le d$, such that the limit
\begin{equation}
	\label{e2.4}
	\frac{f\big(t^EL(t)x\big)}{t^{-\text{tr}(E)}V(t)}\to \lambda(x), 
\end{equation}
holds locally uniformly in $x\in R_\epsilon:=\{x: ||x||\ge \epsilon\}\subset \mathbb{R}^d$, $\epsilon>0$, where $||\cdot||$ is a norm on $\mathbb{R}^d$. Since $E=\mbox{DIAG}(\lambda_i)$, where $\lambda_i>0$, $1\le i\le d$, $t^E = \mbox{DIAG}(t^{\lambda_i})$, the diagonal matrix with diagonal entries $t^{\lambda_i}$s. 

Let $[x] := \sum_{i=1}^d |x_i|^{1/\lambda_i}$, $x\in \mathbb{R}^d_+$. While $[x]$ is not a norm, the function is known as a quasi-homogeneous function, with scaling $[t^Ex] = t[x]$, $t>0$. Since $[x]$ is unbounded if at least one of $x_i$s go to $\pm\infty$, the set $Q=\{x\in {\cal C}: [x]=1\}$ is compact. Therefore, the limit \eqref{e2.4} holds uniformly on $Q\cap R_\epsilon$. 

Let $h(t) = t^{-\text{tr}(E)}V(t)\in \mbox{RV}_{-\text{tr}(E)-\rho}$. Consider any Borel subset $B$ of $R_\epsilon$ that is bounded away from zero, $\epsilon>0$. For $x\in B$, 
\[\frac{f(t^EL(t)x)}{h(t)}= \frac{f(t^E[x]^E\,L(t)\,[x]^{-E}x)}{h(t[x])}\frac{h(t[x])}{h(t)}. 
\]
Since $[x]^{-E}x\in Q$ and \eqref{e2.4} holds uniformly on $Q\cap R_\epsilon$, for any given $\delta>0$, when $t[x]>t_1$, 
\[\frac{f\big((t[x])^E\,L(t)\,[x]^{-E}x\big)}{h(t[x])}\le \sup_{x\in B}\lambda\big([x]^{-E}x\big) + \delta \le \sup_{x\in Q}\lambda(x)+\delta = k <\infty,
\]
which follows from the compactness of $Q$ and local boundedness of $\lambda(x)$. Let $\epsilon_0$ be the smallest value of $[x]$ on $R_\epsilon$, and obviously $\epsilon_0>0$. 
Therefore, whenever $t>t_1\epsilon_0^{-1}$, 
\begin{equation}
	\label{e2.5}
	\sup_{x\in B}\frac{f\big((t[x])^E\,L(t)\,[x]^{-E}x\big)}{h(t[x])}\le \sup_{x\in B}\lambda\big([x]^{-E}x\big) + \delta \le \sup_{x\in Q}\lambda(x)+\delta = k <\infty.
\end{equation}
In addition, since $h(t)\in \mbox{RV}_{-\text{tr}(E)-\rho}$, it follows from Karamata's representation, the uniform convergence theorem, and quasi-homogeneity that whenever $t>t_2$, for all $[x]>1$, 
\begin{equation}
	\label{e2.6}
	\frac{h(t[x])}{h(t)}\le  c[x]^{-\text{tr}(E)-\rho+\gamma}, \ x\in B,
\end{equation}	
for any small $0<\gamma<\rho$ and a constant $c>0$. Therefore, it follows from \eqref{e2.5} and \eqref{e2.6} that, whenever $t>\max\{t_1\epsilon_0^{-1},t_2\}$, $f\big(t^EL(t)x\big)/h(t)\le \kappa_1  [x]^{-\text{tr}(E)-\rho+\gamma}$, $x\in B$ and $[x]>1$, where $\kappa_1>0$ is a constant. 

To show that $[x]^{-\text{tr}(E)-\rho+\gamma}$ is Lebesque integrable on $R_\epsilon$, for any small $\epsilon>0$, consider the following decomposition:
\[\{x: ||x||\ge  \epsilon\} = \mbox{I}+\mbox{II}, \  \ \mbox{I}\cap \mbox{II}=\emptyset, 
\]
where $\mbox{I} = \{x: |x_i|\le 1, 1\le i\le d\}\cap R_\epsilon$ is compact, and $\mbox{II}=\{x: |x_j|>1, j\in A; x_i\le 1, i\in A^c\}\cap R_\epsilon$, for some $\emptyset \ne A\subseteq \{1, \dots, d\}$. Obviously $[x]^{-\text{tr}(E)-\rho+\gamma}$ is bounded in $\mbox{I}$ and thus integrable on $\mbox{I}$. Observe that 
\by
& & \int_{\mbox{II}} \Big(\sum_{i=1}^d |x_i|^{1/\lambda_i}\Big)^{-\text{tr}(E)-\rho+\gamma}dx \\
&\le  & \int_{\{x: |x_j|>1, j\in A\}}\int_{\{x: 0\le |x_i|\le 1, i\in A^c\}}\Big(\sum_{i=1}^d |x_i|^{1/\lambda_i}\Big)^{-\text{tr}(E)-\rho+\gamma}dx_{A^c}dx_A\\
&\le &  2^{|A^c|}\int_{\{x: |x_j|>1, j\in A\}}\Big( \sum_{j\in A}|x_j|^{1/\lambda_j}\Big)^{-\text{tr}(E)-\rho+\gamma}dx_{A}
\ey
\by
&= & 2^{|A^c|}\int_{\{y: |y_j|>1, j\in A\}}\Big(\sum_{j\in A}|y_j|\Big)^{-\text{tr}(E)-\rho+\gamma}\prod_{j\in A}\lambda_j|y_j|^{\lambda_j-1}dy_{A}\\
&\le &  2^{|A^c|}|A|^{-\text{tr}(E)-\rho+\gamma} \int_{\{y: |y_j|>1, j\in A\}}y_{(|A|)}^{-\text{tr}(E)-\rho+\gamma} y_{(|A|)}^{\sum_{j\in A}\lambda_j-|A|}\prod_{j\in A}\lambda_j\,dy_{A}\\
&= &   2^{|A^c|}|A|^{-\text{tr}(E)-\rho+\gamma}\prod_{j\in A}\lambda_j \int_{\{y: |y_j|>1, j\in A\}}\big(y_{(|A|)}\big)^{-|A|-\text{tr}(E)+\sum_{j\in A}\lambda_j-\rho+\gamma}dy_{A}<\infty, 
\ey
which follows from the facts that $y_{(|A|)}=\max\{|y_j|, j\in A\}$ is the $L_\infty$-form on $\mathbb{R}^{|A|}$ and $-|A|-\text{tr}(E)+\sum_{j\in A}\lambda_j-\rho+\gamma<-|A|$. 

Since the function $\kappa_1 [x]^{-\text{tr}(E)-\rho+\gamma}$ is Lebesque integrable on $R_\epsilon$, it then follows from dominated convergence that $\lambda(x)$ is Lebesque integrable on $B\subseteq R_\epsilon$ and 
\[\frac{\mathbb{P}(X\in t^EL(t)B)}{U(t)}=\int_B \frac{f(t^EL(t)x)}{t^{-\text{tr}(E)}V(t)}dx\to \int_B\lambda(x)dx
\]
where $U(t) = V(t)\big(\prod_{i=1}^d\ell_i(t)\big)\in \mbox{RV}_{-\rho}$, and \eqref{e2.3} follows. \hfill $\Box$

\begin{Rem}\label{operator-r-2}
	\begin{enumerate}
		\item Let $\Lambda(B):= \int_B \lambda(x)dx$ for every relative compact set $B\subseteq \overline{\mathbb{R}}^d\backslash \{0\}$. Theorem \ref{pdf->CDF} states that if a density is operator-regularly varying, then its distribution is operator-regularly varying with intensity measure $\Lambda(\cdot)$. The result was obtained in \cite{Li2023} for scaling function $t^E$ and Theorem \ref{pdf->CDF} is proved for general scaling $t^EL(t)$, where $L(t)$ is a diagonal matrix with slowly varying entries. 
		\item Suppose that $X=(X_1, \dots, X_d)$ is non-negative. Let $B_i = [0,\infty]^{i-1}\times (x_i,\infty]\times [0,\infty]^{d-i}$, $x_i>0$, be relative compact in $\overline{\mathbb{R}}^d_+\backslash \{0\}$, $1\le i\le d$, with respect to the topology generated by the Alexandroff compactification-uncompactification. It follows from Theorem \ref{pdf->CDF} that for $1\le i\le d$, 
			\begin{equation}
			\label{operator-e-6}
			\frac{\mathbb{P}(X_i> t^{\lambda_i}\ell_i(t)x_i)}{U(t)}\to \int_{B_i}\lambda(x)dx,
		\end{equation}
	for $\lambda_i>0$ and $\ell_i(t)\in \mbox{RV}_0$. That is, the $i$-th marginal distribution is regularly varying with tail index $-\rho/\lambda_i$. 
	\end{enumerate}
\end{Rem}

For multivariate densities, pointwise convergence implies convergence in total variation by the classical Scheffé theorem, whose proof relies on the Riesz representation (duality) theorem. In particular, Riesz duality identifies the 
$L^1$-norm with the operator norm induced by bounded linear functionals represented by elements of 
$L^\infty$. Scheffé’s argument then exploits a maximizing sequence of indicator functions of measurable sets, together with the dominated convergence theorem, to obtain total variation convergence.

However, the Riesz duality between 
$L^1$
and 
$L^\infty$
does not in general extend to more general spaces. In Theorem~\ref{pdf->CDF}, we instead establish distributional convergence by combining local uniform convergence and a scaling property with the dominated convergence theorem. From this perspective, Theorem~\ref{pdf->CDF} may be viewed as a generalization of Scheffé’s theorem to the setting of Radon measures.

\section{Decomposition of Operator-Regular Variation}

While operator regular variation with diagonal norming matrices captures extremal dependence structures of multivariate distributions without requiring tail equivalence, the associated underlying copula always exhibits a tail density of order $\kappa = (1, \dots,1)$.

\begin{The}
	\label{operator-t-1} Let $F$ be a distribution defined on $\mathbb{R}^d_+$, with ultimately non-increasing density $f$, and let $C$ denote its copula with density $c$. If 
	$f\in \mbox{MRV}(E, \rho,\lambda(\cdot))$, where $E$ is given by \eqref{diag}, $\lambda_i>0$, $1\le i\le d$, then the upper tail density $\lambda_C(\cdot; (1, \dots, 1))$ exists and  
	\begin{equation}\label{operator-e-9}
		 \lambda(w) = \lambda_C\big((w_1^{-\alpha_1}, \dots, w_d^{-\alpha_d}); (1, \dots, 1)\big)|J(w_1^{-\alpha_1}, \dots, w_d^{-\alpha_d})|,
	\end{equation}
for $w_i>0$, $1\le i\le d$, where $|J(w_1^{-\alpha_1}, \dots, w_d^{-\alpha_d})|=\prod_{i=1}^d \alpha_iw_i^{-\alpha_i-1}$ is the Jacobian determinant of the homeomorphic transform $y_i=w_i^{-\alpha_i}$, $\alpha_i = \rho/\lambda_i$, $1\le i\le d$.
\end{The}

\noindent
{\sl Proof.} Assume, without loss of generality, that $\rho = 1$, and $V(t) = t^{-1}\ell(t)$ in \eqref{e2.1}, where $\ell(t)\in \mbox{RV}_0$. The density $f$ of $F$ is given by 
\begin{equation}
	\label{f-c}
	f({x})=c(F_1(x_1), \dots, F_d(x_d))\prod_{i=1}^df_i(x_i),~{ x}=(x_1, \dots, x_d)\in \mathbb{R}^d_+, 
\end{equation}
where $F_i$ is the marginal distributions of $F$ with density $f_i$, $1\le i\le d$, 
and thus 
\[
f(g(t)x)=c\big(F_1(g_1(t)x_1), \dots, F_d(g_d(t)x_d)\big)\prod_{i=1}^df_i\big(g_i(t)x_i\big),~t>0,~{x}=(x_1, \dots, x_d), 
\]
where $g(t)=\mbox{DIAG}(g_i(t))$ with $g_i(t)=t^{\lambda_i}l_i(t)$, $l_i \in \mbox{RV}_0$, $1\le i\le d$. By Theorem \ref{pdf->CDF}, $F_i\in \mbox{RV}_{-\alpha_i}$, where $\alpha_i = 1/\lambda_i$, $1\le i\le d$. Since $f$ is ultimately non-increasing, the marginal densities $f_i$, $1\le i\le d$, are ultimately non-increasing. It then follows from Landau's theorem (see \cite{Resnick07}) that $f_i\in \mbox{RV}_{-\alpha_i-1}$; that is,  $f_i(t)=t^{-\alpha_i-1}\ell_i(t)$, where $\ell_i\in \mbox{RV}_0$, $\alpha_i>0$, $1\le i\le d$.  
It follows from Karamata's Theorem (see, e.g.,  Theorem 2.1 in \cite{Resnick07}) that for $1\le i\le d$, 
\[\Fbar_i(t)\sim \alpha_i^{-1}tf_i(t) \in \RV_{-\alpha_i}, ~\mbox{and}~\Fbar_i(tx_i)/\Fbar_i(t)\sim x_i^{-\alpha_i},
\]
where $\overline{F}_i(t)=1-F_i(t)$ denotes the survival function of $F_i$, $1\le i\le d$. 

Consider 
\by
\frac{f(g(t)w)}{t^{-\text{tr}(E)}V(t)}
&=& \frac{c\left(1-\frac{\Fbar_1(g_1(t)w_1)}{\Fbar_1(g_1(t))}\Fbar_1(g_1(t)),\dots,1-\frac{\Fbar_d(g_d(t)w_d)}{\Fbar_d(g_d(t))}\Fbar_d(g_d(t)) \right)\prod_{i=1}^d f_{i}(g_i(t)w_i)}{t^{-\sum_{i=1}^d\lambda_i}t^{-1}\ell(t) }, 
\ey
where $\ell\in \mbox{RV}_0$. Observe that for some function $s(\cdot)$, 
\by
\frac{c\left(1-\frac{\Fbar_i(g_i(t)w_i)}{\Fbar_i(g_i(t))}\Fbar_i(g_i(t)),\ \forall i \right)}{t^{-1+d}s(t) } &=& \frac{f(g(t)w)t^{-\sum_{i=1}^d\lambda_i}\ell(t)}{t^{d}s(t)t^{-\text{tr}(E)}V(t)\prod_{i=1}^d f_{i}(g_i(t)w_i)}\\
&= & \frac{1}{\prod_{i=1}^d w_i^{-\alpha_i-1}}\frac{f(g(t)w)\ell(t)}{s(t)t^{-\text{tr}(E)}V(t)\prod_{i=1}^d l_i^{-\alpha_i-1}(t)\ell_i(g_i(t)w_i)}. 
\ey
Let $s(t) := \ell(t)\big[\prod_{i=1}^d\alpha_i l_i^{-\alpha_i-1}(t)\ell_i(g_i(t))\big]^{-1}$. Obviously $s\in \mbox{RV}_0$ and 
\[s(t)\sim \ell(t)\big[\prod_{i=1}^d\alpha_i l_i^{-\alpha_i-1}(t)\ell_i(g_i(t)w_i)\big]^{-1}
\]
locally uniformly in $w_i> 0$, $1\le i\le d$. Therefore, 
\[\frac{c\left(1-\frac{\Fbar_i(g_i(t)w_i)}{\Fbar_i(g_i(t))}\Fbar_i(g_i(t)),\ \forall i \right)}{t^{-1+d}s(t) }\sim  \frac{1}{\prod_{i=1}^d \alpha_iw_i^{-\alpha_i-1}}\frac{f(g(t)w)}{t^{-\text{tr}(E)}V(t)}
\]
locally uniformly in $w_i>0$, $1\le i\le d$, which implies that 
\begin{equation}\label{operator-e-11}
\frac{c\left(1-\frac{\Fbar_i(g_i(t)w_i)}{\Fbar_i(g_i(t))}\Fbar_i(g_i(t)),\ \forall i \right)}{t^{-1+d}s(t) }\to \frac{1}{\prod_{i=1}^d \alpha_iw_i^{-\alpha_i-1}}\lambda(w)	
\end{equation}
locally uniformly in $w_i>0$, $1\le i\le d$. 

Since $\frac{\Fbar_i(g_i(t)w_i)}{\Fbar_i(g_i(t))}\to w_i^{-\alpha_i}$, $1\le i\le d$, locally uniformly, \eqref{operator-e-11} implies, via the sequential property of uniformity, that 
\begin{equation}\label{operator-e-12}
	\frac{c\left(1-\frac{\Fbar_i(g_i(t)w_i)}{\Fbar_i(g_i(t))}\Fbar_i(g_i(t)),\ \forall i \right)}{t^{-1+d}s(t) }\sim \frac{c\left(1-w_i^{-\alpha_i}\Fbar_i(g_i(t)),\ \forall i \right)}{t^{-1+d}s(t) }
\end{equation}
locally uniformly in $w_i>0$, $1\le i\le d$. Putting \eqref{operator-e-11} and \eqref{operator-e-12} together, we have
\begin{equation}
\label{operator-e-13}	
\frac{c\left(1-w_i^{-\alpha_i}\Fbar_i(g_i(t)),\ \forall i \right)}{t^{-1+d}s(t) }\to \frac{1}{\prod_{i=1}^d \alpha_iw_i^{-\alpha_i-1}}\lambda(w)	
\end{equation}
locally uniformly in $w_i>0$, $1\le i\le d$. Let $u=t^{-1}$, and $u\to 0$ if and only if $t\to \infty$. Rewrite \eqref{operator-e-13}, leading to
\[\lim_{u\to 0}\frac{c\left(1-w_i^{-\alpha_i}\Fbar_i(g_i(u^{-1})),\ \forall i \right)}{u^{1-d}s(u^{-1}) }= \frac{1}{\prod_{i=1}^d \alpha_iw_i^{-\alpha_i-1}}\lambda(w)	
\]
locally uniformly in $w_i>0$, where $s(u^{-1})\in \mbox{RV}_0(0)$ and $\Fbar_i(g_i(u^{-1}))\in \mbox{RV}_{1}(0)$, $1\le i\le d$. This shows that the upper tail density $\lambda_C(\cdot; (1, \dots, 1))$ exists and \eqref{operator-e-9} holds. 
 \hfill $\Box$

\begin{Rem}\label{operator-r-5}
	\begin{enumerate}
		\item 	Since $f\in \mbox{MRV}(E, -\rho,\lambda(\cdot))$, $\rho>0$, if and only if $f\in \mbox{MRV}(E/\rho, -1,\lambda(\cdot))$, then if 
		$f\in \mbox{MRV}(E, \rho,\lambda(\cdot))$, where $E$ is given by \eqref{diag}, $\lambda_i>0$, $1\le i\le d$, then the upper tail density is given by \eqref{operator-e-9}, where $\alpha_i=\rho/\lambda_i$, $1\le i\le d$. 
		\item If the copula density $c(\cdot)$ is continuous in a small open neighborhood of  $(1, \dots, 1)$, and \eqref{operator-e-1} holds locally uniformly, the upper tail density $\lambda_C(\cdot;\kappa)$ is continuous. By continuous extensions, $\lambda_C(\cdot;\kappa)$ is well-defined on the entire $\mathbb{R}_+^d$. 
		\item The tail order vector 
		$\kappa = (1, \dots, 1)$ reflects the fact that multivariate regular variation with a common operator norming extends over the entire cone 
	 $\overline{\mathbb{R}}_+^d\backslash\{0\}$, including the marginal, univariate sub-cones. In contrast, other tail order vectors arise when multivariate regular variation, under suitable operator norming, is confined to a proper multivariate sub-cone of 
 $\overline{\mathbb{R}}_+^d\backslash\{0\}$.
	\end{enumerate}
\end{Rem}

On the other hand, one can construct an operator regular variation from a copula with tail densities together with univariate regularly varying margins. However, this construction requires a compatibility condition between the copula’s tail densities and the marginal tail behavior. A copula $C$ with upper tail density \eqref{operator-e-1} of order $(\rho_1, \dots, \rho_d)$ is said to be regularly varying compatible with univariate margins $F_i(x)$, with regularly varying tail index $\alpha_i$, $1\le i\le d$, if the scaling functions $r_i(t^{-1})\sim 1-F_i(t^{\rho_i/\alpha_i})$, as $t\to \infty$, $1\le i\le d$.

\begin{The}
	\label{operator-t-2} Let a copula $C$ with density $c(\cdot)$ have the upper tail density $\lambda_C(\cdot; \kappa)$, $\kappa=(\rho_1, \dots, \rho_d)$, and $f_i\in \mbox{RV}_{-\alpha_i-1}$ are the univariate densities, $\rho_i>0$, $\alpha_i>0$, $1\le i\le d$. 
	If $F(x_1, \dots, x_d) = C(F_1(x_1), \dots, F_d(x_d))$, where $C$ and $F_i$, $1\le i\le d$, are compatible, then $F$ is operator-regularly varying with tail density $\lambda(\cdot)$, as defined in \eqref{e2.1}, satisfying that  for any ${w}=(w_1, \dots, w_d)>{0}$, 
	\begin{eqnarray}
		\lambda({w}):=\lim_{t\to \infty}\frac{f(g(t){w})}{t^{-\sum_{i=1}^d\lambda_i}V(t)}
		&=&\lambda_C(w_1^{-\alpha_1}, \dots, w_d^{-\alpha_d}; \kappa)|J(w_1^{-\alpha_1}, \dots, w_d^{-\alpha_d})|,\label{tail density transform}
	\end{eqnarray}
	where $\lambda_i=\rho_i/\alpha_i>0$, $1\le i\le d$,  $g: \mathbb{R}_+\to \mbox{GL}(\mathbb{R}^d)$ is operator-regularly varying with index $E=\mbox{DIAG}(\lambda_i)$ given by \eqref{diag}, $V\in \mbox{RV}_{-1}$ and $J(w_1^{-\alpha_1}, \dots, w_d^{-\alpha_d})=\prod_{i=1}^d\alpha_iw_i^{-\alpha_i-1}$ is the Jacobian determinant of the homeomorphic transform $y_i=w_i^{-\alpha_i}$, $1\le i\le d$.
\end{The}

\noindent
{\sl Proof.}
The density $f$ of $F$ is given by 
\[
	f({x})=c(F_1(x_1), \dots, F_d(x_d))\prod_{i=1}^df_i(x_i),~{ x}=(x_1, \dots, x_d)\in \mathbb{R}^d_+, 
\]
where $F_i$ is the marginal distributions of $F$ with density $f_i$, $1\le i\le d$, 
and thus 
\[
f(g(t)x)=c\big(F_1(g_1(t)x_1), \dots, F_d(g_d(t)x_d)\big)\prod_{i=1}^df_i\big(g_i(t)x_i\big),~t>0,~{x}=(x_1, \dots, x_d), 
\]
where $g(t)=\mbox{DIAG}(g_i(t))=t^E$ with $g_i(t)=t^{\lambda_i}$, $\lambda_i=\rho_i/\alpha_i>0$, $1\le i\le d$. It follows from Karamata's Theorem (see, e.g.,  Theorem 2.1 in \cite{Resnick07}) that for $1\le i\le d$, 
\[\Fbar_i(t)\sim \alpha_i^{-1}tf_i(t) \in \RV_{-\alpha_i}, ~\mbox{and}~\Fbar_i(tx_i)/\Fbar_i(t)\sim x_i^{-\alpha_i}, 
\] 
where $f_i(t) = t^{-\alpha_i-1}\ell_i(t)$, $\ell_i\in \mbox{RV}_0$, $1\le i\le d$. 

Let $u=t^{-1}$, and due to the compatibility, for $1\le i\le d$, $\Fbar_i(g_i(t)) \sim r_i(u)$, as $u\to 0$. Therefore, 
\begin{equation}
	\label{operator-e-14}
\lambda_{C}(w; \kappa)=\lim_{u\to 0}\frac{c\big(1-\Fbar_1(g_1(u^{-1}))w_1, \dots, 1-\Fbar_d(g_d(u^{-1}))w_d\big)}{u^{1-\sum_{i=1}^d\rho_i}\ell(u)}, 
\end{equation}
and the convergence holds locally uniformly in $w \in [0,\infty)^d\backslash\{0\}$. 
Consider, for some function $V(\cdot)$, 
\by
\frac{f(g(t)w)}{t^{-\sum_{i=1}^d\lambda_i}V(t)}
&=& \frac{c\left(1-\frac{\Fbar_1(g_1(t)w_1)}{\Fbar_1(g_1(t))}\Fbar_1(g_1(t)),\dots,1-\frac{\Fbar_d(g_d(t)w_d)}{\Fbar_d(g_d(t))}\Fbar_d(g_d(t)) \right)\prod_{i=1}^d f_{i}(g_i(t)w_i)}{t^{-\sum_{i=1}^d\lambda_i}V(t)}
\ey
\begin{equation}\label{operator-e-15}
=\frac{c\left(1-\frac{\Fbar_1(g_1(t)w_1)}{\Fbar_1(g_1(t))}\Fbar_1(g_1(t)),\dots,1-\frac{\Fbar_d(g_d(t)w_d)}{\Fbar_d(g_d(t))}\Fbar_d(g_d(t)) \right)\prod_{i=1}^d f_{i}(g_i(t)w_i)}{t^{-1+\sum_{i=1}^d\rho_i}\ell(t^{-1})t^{-\sum_{i=1}^d\lambda_i}V(t)t^{1-\sum_{i=1}^d\rho_i}[\ell(t^{-1})]^{-1}}.
\end{equation} 
It follows from the sequential property of local uniformity \eqref{operator-e-14} that 
\be
\frac{c\left(1-\frac{\Fbar_i(g_i(t)w_i)}{\Fbar_i(g_i(t))}\Fbar_i(g_i(t)),\ \forall\ i \right)}{t^{-1+\sum_{i=1}^d\rho_i}\ell(t^{-1})}&\sim &  \frac{c\left(1-w_1^{-\alpha_1}\Fbar_1(g_1(t)),\dots,1-w_d^{-\alpha_d}\Fbar_d(g_d(t)) \right)}{t^{-1+\sum_{i=1}^d\rho_i}\ell(t^{-1})}\nonumber \\
&\to & \lambda_{C}((w_1^{-\alpha_1}, \dots, w_d^{-\alpha_d}); \kappa)\label{operator-e-16}
\ee
locally uniformly in $w_i>0$, $1\le i\le d$, where $\kappa = (\rho_1, \dots, \rho_d)$. On the other hand, 
\by
\frac{\prod_{i=1}^d f_{i}(g_i(t)w_i)}{t^{-\sum_{i=1}^d\lambda_i}V(t)t^{1-\sum_{i=1}^d\rho_i}[\ell(t^{-1})]^{-1}}&\sim & \frac{\big(\prod_{i=1}^dw_i^{-\alpha_i-1}\big)t^{-\sum_{i=1}^d\rho_i-\sum_{i=1}^d\lambda_i}\prod_{i=1}^d\ell_i(t^{\lambda_i}w_i)}{t^{-\sum_{i=1}^d\lambda_i}V(t)t^{1-\sum_{i=1}^d\rho_i}[\ell(t^{-1})]^{-1}}\\
&\sim & \frac{\big(\prod_{i=1}^d\alpha_iw_i^{-\alpha_i-1}\big)\prod_{i=1}^d\ell_i(t^{\lambda_i}w_i)}{V(t)t[\ell(t^{-1})]^{-1}\prod_{i=1}^d\alpha_i}.
\ey
Let $V(t) := t^{-1}\ell(t^{-1})\prod_{i=1}^d\alpha_i^{-1}\ell_i(t^{\lambda_i})$, and clearly $V\in \mbox{RV}_{-1}$, and 
\begin{equation}
\label{operator-e-17}	
\frac{\prod_{i=1}^d f_{i}(g_i(t)w_i)}{t^{-\sum_{i=1}^d\lambda_i}V(t)t^{1-\sum_{i=1}^d\rho_i}[\ell(t^{-1})]^{-1}}\to \prod_{i=1}^d\alpha_iw_i^{-\alpha_i-1}
\end{equation}
locally uniformly in $w_i>0$, $1\le i\le d$, due to local uniformity of univeriate regular variations. Plugging \eqref{operator-e-16} and \eqref{operator-e-17} into \eqref{operator-e-15} yields the result. 
 \hfill $\Box$

\begin{Rem}\label{operator-r-8}
	\begin{enumerate}
		\item 	If a copula $C$ with upper tail density \eqref{operator-e-1} is not regularly varying compatible with univariate margins $F_i(x)$, $1\le i\le d$, then the construction $F(x_1, \dots, x_d)=C(F_1(x_1), \dots, F_d(x_d))$ yields multivariate extremes that may not be amenable to operator regular variation, and can be analyzed only through the underlying copula and univariate regular variations.
				\item If the copula density $c(\cdot)$ is continuous in a small open neighborhood of  $(1, \dots, 1)$, and \eqref{operator-e-1} holds locally uniformly, the upper tail density $\lambda_C(\cdot;\kappa)$ is continuous. It follows from \eqref{tail density transform} that $\lambda(w)$ is continuous in $w>0$. By continuous extensions, the tail density $\lambda(\cdot)$ is well-defined on the cone $\overline{\mathbb{R}}_+^d\backslash\{0\}$. 
	\end{enumerate}
\end{Rem}

The tail dependence functions of a copula $C$ were introduced in \cite{Li2009, JLN10} to describe the scale-invariant tail dependence information of $X=(X_1, \dots, X_d)$ with copula $C$ and univariate marginal distributions $F_i$, $1\le i\le d$. The equivalence between the intensity measures of $X$ with regularly varying distribution $F$ and the tail dependence functions of its underlying copula $X$ were discussed in \cite{LS2009, HJL12, LH13}. The upper operator exponent function with tail order vector $\kappa=(\rho_1, \dots, \rho_d)$, $\rho_i>0$, $1\le i\le d$, is defined as the follows, for every $w=(w_1, \dots, w_d)\in \mathbb{R}_+^d\backslash\{0\}$, 
\begin{equation}\label{operator-e25}
	a_C((w_1, \dots, w_d); \kappa) = \lim_{u\to 0^+}\frac{\mathbb{P}\big(F_i(X_i)> 1-r_i(u)w_i, \mbox{for some}\ 1\le i\le d\big)}{u\ell(u)}
\end{equation} 
for some functions $r_i\in \mbox{RV}_{\rho_i}(0)$, $1\le i\le d$, and $\ell\in \mbox{RV}_0(0)$. The lower operator exponent function of a copula is defined as the upper operator exponent function of its survival copula.  See \cite{Li2018} for details on operator tail dependencies of copulas and their relations with operator regular variation.

\begin{The}
	\label{operator-t-3}
Let a copula $C$ with continuous density $c(\cdot)$ have the upper exponent function $a_C(\cdot; \rho)$ and upper tail density $\lambda_C(\cdot; \rho)$, where $\rho=(\rho_1, \dots, \rho_d)$, $\rho_i>0$, $1\le i\le d$. Then  
\begin{equation}\label{operator-e-19}
	a_C((w_1, \dots, w_d); \rho) = \int_{[0,w]^c}\lambda_C(x;\rho)dx
\end{equation}
for any $w=(w_1, \dots, w_d)\in \mathbb{R}_+^d\backslash\{0\}$.
\end{The}

\noindent
{\sl Proof.}
Let $F_i$ denote a univariate distribution with density $f_i$, $1\le i\le d$. Construct $F_i$ and $f_i$, that are compatible with $C$, such that $F_i\in \mbox{RV}_{-\rho_i}$ and $f_i\in \mbox{RV}_{-\rho_i-1}$, $1\le i\le d$. Therefore, $F(x_1, \dots, x_d)=C(F_1(x_1), \dots, F_d(x_d))$ is a distribution on $\mathbb{R}_+^d$, that is operator-regularly varying in the sense of \eqref{e2.1}, where $E$ is the identity matrix. By Theorem \ref{operator-t-2}, the tail density of $F$ is given by 
\[\lambda({w})
=\lambda_C(w_1^{-\rho_1}, \dots, w_d^{-\rho_d}; \rho)|J(w_1^{-\rho_1}, \dots, w_d^{-\rho_d})|,
\]
with continuous extensions to all $w=(w_1, \dots, w_d)\in \mathbb{R}_+^d\backslash\{0\}$. 
It follows from \cite{Li2018} and Theorem \ref{pdf->CDF} that 
\by
a_C((w_1^{-\rho_1}, \dots, w_d^{-\rho_d}); \rho) &=& \int_{[0,w]^c}\lambda_C(x_1^{-\rho_1}, \dots, x_d^{-\rho_d}; \rho)|J(x_1^{-\rho_1}, \dots, x_d^{-\rho_d})|dx\\
&=& \int_{[0,w]^c}\lambda_C(x_1^{-\rho_1}, \dots, x_d^{-\rho_d}; \rho)dx_1^{-\rho_1} \cdots dx_d^{-\rho_d}\\
&=& \int_{[0,w^{-\rho}]^c}\lambda_C(x_1, \dots, x_d; \rho)dx_1 \cdots dx_d
\ey
where $[0,w^{-\rho}] = \prod_{i=1}^d[0, w_i^{-\rho_i}]$, 
leading to \eqref{operator-e-19}. 
  \hfill $\Box$

\begin{Rem}
\label{operator-r-9}
\begin{enumerate}
	\item The integral relation \eqref{operator-e-19} holds on the entire cone $\mathbb{R}_+^d\backslash\{0\}$. Since \eqref{e2.3} holds for a Borel subset in any sub-cone of $\mathbb{R}_+^d\backslash\{0\}$ and the relations \eqref{operator-e-9} and \eqref{tail density transform} are local, the integral relation between the tail dependence measure generated by $a_C(\cdot; \rho)$ and the tail density holds in any sub-cone as well. 
	\item Assuming the differentiablity, \eqref{operator-e-19} is equivalent to 
	\[\frac{\partial^d a_C((w_1, \dots, w_d); \rho)}{\partial w_1\cdots \partial w_d} = \lambda_C(w;\rho)
	\]
	for any $w=(w_1, \dots, w_d)\in \mathbb{R}_+^d\backslash\{0\}$. Such a differential form was used for higher order tail dependence functions in \cite{HJL12,LH14}. 
\end{enumerate}
\end{Rem}

\section{An Example of Liouville Copulas with Operator Tail Densities}
\label{Liouville Copulas}

An absolutely continuous non-negative random vector $X=(X_1, \dots, X_d)$ with distribution $F$ is said to
have a Liouville distribution, denoted by $X\sim L_d[g(t); a_1, \dots, a_d]$ or $F\sim L_d[g(t); a_1, \dots, a_d]$,  if the its joint probability
density function is proportional to
\begin{equation}
	\label{Liouville0}
	g\Big(\sum_{i=1}^dx_i\Big)\prod_{i=1}^dx_i^{a_i-1}
\end{equation}
for $x_1>0, \dots, x_d>0$, where $a_i>0$, $i=1, \dots, d$, and the driving function $g(\cdot)$ is a suitably chosen non-negative continuous function, satisfying the integrablibity that
\begin{equation}
	\label{Liouville1}
	\int_0^\infty t^{\sum_{i=1}^da_i-1}g(t)dt < \infty. 
\end{equation}
This condition is assumed to ensure that \eqref{Liouville0} is a probability density function, due to the well-known formula for Liouville's integral. 
We also assume in this section  that $g(\cdot)$ has the non-compact support $[0, \infty)$. 
For example, the inverted Dirichlet distribution has the joint density function 
\[f(x_1, \dots, x_d) = \frac{\Gamma\big(\sum_{i=1}^{d+1}a_i\big)}{\Gamma(a_{d+1})}\Big(1+\sum_{i=1}^dx_i\Big)^{-a_1- \cdots -a_d-a_{d+1}}\prod_{i=1}^d\frac{x_i^{a_i-1}}{\Gamma(a_i)}, 
\]
in which, $g(t) = (1+t)^{-a_1-\dots -a_d-a_{d+1}}$, $t>0$, $a_i>0$, $i=1, \dots, d+1$, where $\Gamma(\cdot)$ denotes the gamma function. Observe that this function is regularly varying with tail index $\sum_{i=1}^{d+1}a_i$. In general, however, the function $g(\cdot)$ can be any non-negative function, including rapidly varying functions.  

In general, the process of conditioning a random vector on the sum of its components leading a distribution of certain independent events can be modeled using the multivariate Liouville distribution and its various extensions. The theory and extensive discussions of multivariate Liouville distributions can be found in \cite{Gupta87, Gupta91, Gupta92}, and the history and related references are detailed in \cite{Gupta01}. 

If $X=(X_1, \dots, X_d)\sim L_d[g(t); a_1, \dots, a_d]$ with the joint distribution $F(x)$ on $\mathbb{R}_+^d$, then $X_1\sim L_1[g_1(t); a_1]$ with marginal distribution $F_1$, where $g_1(\cdot)$ the Weyl fractional
integral of order $a^{(1)} = \sum_{i=2}^da_i$, given by 
\[g_1(t) = W^{a^{(1)}}g(t) = \frac{1}{\Gamma (a^{(1)})}\int_t^\infty (s-t)^{a^{(1)} - 1}g(s)ds,\ t>0.
\]
That is, the margin $X_1$ has the marginal density$f_1(x) = \kappa_1\, W^{a^{(1)}}g(x)\,x^{a_1-1}$, $x>0$, where $\kappa_1>0$ is a normalizing constant. The other marginal densities, however, are not explicit \cite{Gupta87}. A copula $C$ is called a Liouville Copula if $C(u_1, \dots, u_d)= F(F_1^{-1}(u_1), \dots, F_d^{-1}(u_d)))$ for a multivariate Liouville distribution $F\sim L_d[g(t); a_1, \dots, a_d]$, with marginal distributions $F_i$, $1\le i\le d$. A Liouville copula describes the scale-invariant tail dependence of random variables conditioning on their sum, even through it is not explicit.

According to Theorem 3.2 in \cite{Li2023}, if $g\in \mbox{RV}_{-\beta}$, then the density of $X$ is regularly varying with operator tail index $E=\mbox{DIAG}\,(\lambda_i)$, for any $\lambda_i>0$, $i=1, \dots, d$, and with limiting density
\begin{equation}\label{operator-e-18}
	\frac{f(t^{\lambda_1}x_1, \dots, t^{\lambda_d}x_d)}{t^{-\sum_{i=1}^d\lambda_i}V(t)}\to \lambda(x) :=  \Big(\sum_{i\in (\lambda)} x_i\Big)^{-\beta}\prod_{i=1}^dx_i^{a_i-1},
\end{equation}
where $V(t) = g(t^\lambda)t^{\sum_{i=1}^d\lambda_ia_i}\in \mbox{RV}_{-\lambda\beta+\sum_{i=1}^d\lambda_ia_i}$, $\lambda=\max_{1\le k\le d}\{\lambda_k\}$ and $(\lambda) = \{i: \lambda_i = \max_{1\le k\le d}\{\lambda_k\}\}$. It follows from Theorem \ref{operator-t-1} that  a Liouville Copula has an explicit upper tail density as follows
\be
\lambda_C((w_1, \dots, w_d); (1, \dots, 1))&=& \lambda(w_1^{-\frac{1}{\alpha_1}}, \dots, w_d^{-\frac{1}{\alpha_d}})\prod_{i=1}^d\alpha_i^{-1}\prod_{i=1}^d w_i^{-\frac{\alpha_i+1}{\alpha_i}}\nonumber\\
&=& \Big(\sum_{i\in (\lambda)} w_i^{-\frac{1}{\alpha_i}}\Big)^{-\beta}\prod_{i=1}^d\alpha_i^{-1}\prod_{i=1}^dw_i^{-\frac{a_i-1}{\alpha_i}}\prod_{i=1}^d w_i^{-\frac{\alpha_i+1}{\alpha_i}}\nonumber\\
&=&  \Big(\sum_{i\in (\lambda)} w_i^{-\frac{1}{\alpha_i}}\Big)^{-\beta}\prod_{i=1}^d\alpha_i^{-1}\prod_{i=1}^d w_i^{-\frac{a_i-1}{\alpha_i}-\frac{\alpha_i+1}{\alpha_i}}\nonumber\\
&=& \Big(\sum_{i\in (\lambda)} w_i^{-\frac{1}{\alpha_i}}\Big)^{-\beta}\prod_{i=1}^d\alpha_i^{-1}\prod_{i=1}^d w_i^{-\frac{\alpha_i+a_i}{\alpha_i}},\label{operator-e-21}
\ee
where $\alpha_i= \frac{-\lambda\beta+\sum_{i=1}^d\lambda_ia_i}{\lambda_i}$, $1\le i\le d$. 

For a symmetrical Liouville distribution $F\sim L_d[g(t); a, \dots, a]$, where $g\in \mbox{RV}_{-a}$, $a>0$, the upper tail density of its copula $C$ has a simple expression. Since the density of $F$ is regularly varying with operator tail index $E=\mbox{DIAG}\,(a)$, then 
\[\lambda_C((w_1, \dots, w_d); (1, \dots, 1)) = a^{-d}\Big(\sum_{i=1}^d w_i^{-\frac{1}{a}}\Big)^{-a}\prod_{i=1}^dw_i^{-d},
\]
for $w_1>0, \dots, w_d>0$. 

\begin{Rem}
	\begin{enumerate}
	\item Like many copulas, including elliptical copulas, a Liouville copula does not admit a closed-form expression. Its tail density, however, does have the explicit form given in \eqref{operator-e-21}.
	\item A Liouville distribution is not tail-equivalent, but the regular variation in \eqref{operator-e-18} employs operator norming to accommodate different tail decay rates across its margins. Note, however, that the operator tail density of its copula, given in \eqref{operator-e-21}, remains a tail-dependence density of order 1. It is conjectured that the operator tail density of a copula with a diagonal norming matrix whose entries are not all equal may arise in the context of hidden multivariate regular variation within a sub-cone of $\mathbb{R}_+^d\backslash\{0\}$.
\end{enumerate}
\end{Rem}

\end{document}